\documentclass[journal]{IEEEtran}
\usepackage[latin1,utf8]{inputenx}
\usepackage[dvips]{graphicx}
\usepackage{subfigure}
\usepackage[cmex10]{amsmath}
\usepackage{array}
\usepackage[colorlinks=true,citecolor=blue,linkcolor=blue,urlcolor=blue,plainpages=false]{hyperref}
\usepackage{booktabs, tabularx, multirow, array}
\usepackage{nomencl}
\usepackage{multirow}
\usepackage{comment} 
\usepackage{amsmath}
\usepackage{amsfonts}
\usepackage{amssymb}
\usepackage{graphicx}
\usepackage{amsthm}
\usepackage{eurosym}
\usepackage{algorithm}
\usepackage{tabularx}
\usepackage{cite}
\usepackage{array}
\usepackage{graphicx}
\usepackage[cmex10]{amsmath}
\usepackage{amsfonts,amssymb,amsxtra,balance}
\usepackage{layouts}

\usepackage{enumitem}
 \usepackage[table,xcdraw]{xcolor}

\usepackage{epstopdf}
\usepackage{MnSymbol,wasysym}
\usepackage{algorithm}
\usepackage{algpseudocode}
\usepackage{enumitem}

\usepackage{multirow}
\usepackage{dcolumn}
\usepackage{color}

\usepackage{setspace}
\usepackage[export]{adjustbox}
\usepackage{adjustbox}
\usepackage{soul}
\usepackage{booktabs}

\usepackage{tikz}
\usetikzlibrary{shapes,arrows}

\usepackage[textsize=scriptsize,shadow]{todonotes}

\DeclareGraphicsExtensions{.pdf,.png,.jpg}
\graphicspath{{images/}}

\usepackage{marginnote}
\makenomenclature

\begin{document}
\title{An Exact and Scalable Problem Decomposition for Security-Constrained Optimal Power Flow}

\author{Alexandre~Velloso,
	Pascal~Van~Hentenryck,
	Emma~S.~Johnson
}

\maketitle

\begin{abstract}
In this paper, we present decomposition techniques for solving large-scale instances of the security-constrained optimal power flow (SCOPF) problem with primary response. Specifically, under each contingency state, we require that the nodal demands are met and that the synchronized units generating below their limits follow a linear model for primary response. The resulting formulation is a mixed-integer linear program since the primary response model introduces disjunctions to the SCOPF problem. Unfortunately, exact methods relying on traditional Benders’ decomposition do not scale well. As an alternative, we propose a decomposition scheme based on the column-and-constraint-generation algorithm where we iteratively add disjunctions and cuts. We provide procedures for preprocessing dedicated cuts and for numerically determining the post-contingency responses based on the master problem solutions. We also discuss heuristics to generate high-quality primal solutions and upper bounds for the method. Finally, we demonstrate the efficiency of the proposed method on large-scale systems.
\end{abstract}

\begin{IEEEkeywords}
Column and constraint generation, decomposition methods, primary response, security-constrained optimal power flow.
\end{IEEEkeywords}


\section*{Nomenclature}

This section introduces our notation. We use bold symbols for matrices (uppercase) and vectors (lowercase). Additional symbols can be interpreted by the following general rules: Symbols with superscript ``$(j)$" denote new variables, parameters or sets corresponding to the $j$-th iteration of the solution method. The symbols with superscript ``$(*)$" denote the optimal value of the associated (iterating) variable.
	\vspace{0.18cm}

\noindent {\bf Sets}

\begin{description} [labelindent=5.8pt ,labelwidth=35pt, labelsep=4pt, leftmargin =44.8pt, style =standard, itemindent=0pt]
	
	\vspace{0.02cm}

	\vspace{0.033cm}  \item[$\mathcal{E},\mathcal{E}_s$] Feasibility sets for the nominal power flow constraints and for the power flow constraints under contingency state $s$, respectively.
	
	
	\vspace{0.033cm}  \item[$\mathcal{F}_s$] Feasibility set for primary response constraints under contingency state $s$.

	\vspace{0.033cm}  \item[$\mathcal G, \mathcal L, \mathcal N$] Sets of generators, transmission lines and buses, respectively. 
	
	\vspace{0.033cm}  \item[$\mathcal H$] Subset of $\mathcal G$ for devising primal solutions.
		
	
	
	\vspace{0.033cm}  \item[$\mathcal S$] Set of contingencies.
	
	\vspace{0.033cm}  \item[$\mathbb{S}$] Subset of $\mathcal S$ with disjunctive constraints, used in the column-and-constraint-generation algorithm.
	
	\vspace{0.033cm}  \item[$\mathcal{X}, \mathcal{X}_s$] Sets of power flow-related decision variables for nominal state and for contingency state $s$.

	
	\vspace{0.033cm}  \item[$\mathcal{Y}_s$] Set of decision variables associated with primary response under contingency state $s$.
		
	

	
	
	
\end{description}

\vspace{0.9cm}

\noindent {\bf Parameters}

\begin{description} [labelindent=5.8pt ,labelwidth=35pt, labelsep=4pt, leftmargin =44.8pt, style =standard, itemindent=0pt]

	\vspace{0.02cm}  \item[$\alpha$, $\alpha_{s,l}$] Largest transmission line capacity violation and violation for transmission line $l$, under contingency state $s$.

	
	\vspace{0.033cm}  \item[$\beta, \beta_1,\beta_2$] Parameters for selecting preprocessed cuts.
		
	\vspace{0.033cm} \item[$\boldsymbol{\gamma}$] Vector of parameters for primary response.

	\vspace{0.033cm}  \item[$\gamma_i$] Parameter for primary response of generator $i$.

	\vspace{0.033cm}  \item[$\epsilon$] Tolerance for transmission line violation.

	\vspace{0.033cm}  \item[$\varepsilon$] Tolerance for the binary search procedure.	
	
	\vspace{0.033cm}  \item[$\mathbf{A}$, $\mathbf{B}$] Line-bus and Generator-bus incidence matrices.
	

	
	\vspace{0.033cm}  \item[$\mathbf{c}$] Vector of generation costs.
	
	\vspace{0.033cm}  \item[$c_i$] Generation cost of generator $i$.
		
	\vspace{0.033cm}  \item[$\mathbf{{d}}$] Vector of nodal net loads.
	
	\vspace{0.033cm}  \item[$\mathbf{e}$] Vector of ones with appropriate dimension.
	
	\vspace{0.033cm}  \item[$e_s$] Total load imbalance for contingency state $s$.
	
	\vspace{0.033cm}  \item[$\mathbf{\overline{f}}$] Vector of line capacities.
	
	
	\vspace{0.033cm}  \item[$\mathbf{\underline{g}}, \mathbf{\overline{g}}$] Vectors of lower and upper limits for generators.
	
	\vspace{0.033cm}   \item[$\mathbf{\hat{g}}$] Vector of capacities for generators.
	
	\vspace{0.033cm}  \item[$\overline{g}_i$] Upper limit for  generator $i$.
	
	\vspace{0.033cm}  \item[$\hat{g}_i$] Capacity of generator $i$.
	
	\vspace{0.033cm}  \item[$\mathbf{K}_0$] Matrix of power transfer distribution factors.
	
	\vspace{0.033cm}  \item[$\mathbf{K}_1,\mathbf{k}_2$] Preprocessed structures for positive flow limits.
	
	\vspace{0.033cm}  \item[$\mathbf{K}_3,\mathbf{k}_4$] Preprocessed structures for negative flow limits.
	
	\vspace{0.033cm}  \item[$M$] Big-M.
		
	\vspace{0.033cm}  \item[$lb$, $ub$] Lower/upper bound for the decomposition method.
	
	\vspace{0.033cm}  \item[$p$] Parameter for primal solution approach.

	\vspace{0.033cm}  \item[$\mathbf{\overline{r}}$] Vector of primary response limits of generators.
	
	\vspace{0.033cm}  \item[$\overline{r}_i$] Primary response limit of generator $i$, given by $\gamma_i\hat{g}_i$.
	
	\vspace{0.033cm}  \item[$\mathbf{S}$] Angle-to-flow matrix.
	
	\vspace{0.033cm}  \item[$\mathbf{S}^{'}$] Last $|\mathcal{N}|-1$ columns of $\mathbf{S}$.

	\vspace{0.033cm}  \item[$z$] Objective value of the master problem for the column-and-constraint-generation algorithm.
	
	\vspace{0.033cm}  \item[$z_p$] Objective value $z$ when using the parameter $p$.

	
	
	
	
	
	
	
	
	
\end{description}

\vspace{0.2cm}

\noindent {\bf Nominal-state-related decision variables and vectors}

\begin{description} [labelindent=5.8pt ,labelwidth=35pt, labelsep=4pt, leftmargin =44.8pt, style =standard, itemindent=0pt]

	\vspace{0.02cm} \item[$\mathbf{\boldsymbol{\theta}},\mathbf{f},\mathbf{g}$] Phase angles, line flows, and nominal generation.



	\vspace{0.033cm}  \item[$g_i$] Generation of generator $i$ in nominal state.



\end{description}

\vspace{0.08cm}

\noindent {\bf Contingency-state-related decision variables and vectors}

\begin{description}[labelindent=5.8pt ,labelwidth=31pt, labelsep=4pt, leftmargin =40.8pt, style =standard, itemindent=0pt]

	\item[$\mathbf{\boldsymbol{\theta}}_s$] Vector of phase angles under contingency state $s$.

		\vspace{0.033cm}\item[$\mathbf{\boldsymbol{\theta}}_s^{'}$] Last $|\mathcal{N}|-1$ rows of $\mathbf{\boldsymbol{\theta}}_s$.

	\vspace{0.033cm}  \item[$\boldsymbol{\mu}_s$] Vector of dual variables associated with nodal load balance constraint under contingency state $s$.

	\vspace{0.033cm}  \item[$\mathbf{f}_{s}$] Vector for line flows under contingency state $s$.

	\vspace{0.033cm}  \item[$\mathbf{g}_{s}$] Vector for generation under contingency state $s$.
	
	\vspace{0.033cm}  \item[$\mathbf{g}_{s}^{'}$] Provisional vector for $\mathbf{g}_{s}$.
	
	\vspace{0.033cm}  \item[$g_{s,i}$] Generation of generator $i$ under contingency state $s$.
	
	\vspace{0.033cm}  \item[$g_{s,i}^{'}$] Provisional variable for $g_{s,i}$.	
	
	\vspace{0.033cm}  \item[$n_s$] Global signal under contingency state $s$.


	\vspace{0.033cm}  \item[$\mathbf{u}_s^+, \mathbf{u}_s^-$] Vectors of slack variables for line capacities.
	
	\vspace{0.033cm}  \item[$\mathbf{v}_s^+, \mathbf{v}_s^-$] Vectors of slack variables for nodal demand balance.
	
	\vspace{0.033cm}  \item[$\mathbf{x}_{s}$] Binary vector indicating whether generators reached $\mathbf{\overline{g}}$ under contingency state $s$.
	
	\vspace{0.033cm}  \item[$x_{s,i}$] Binary variable indicating whether generator $i$ reached $\overline{g}_i$ under contingency state $s$.


\end{description}

\section{Introduction}

System reliability under contingencies has been widely discussed in the literature. In this context, the goal of the well-known security constrained optimal power flow (SCOPF) problem \cite{ alsac1974optimal ,  bouffard2005umbrella, li2008decomposed,  capitanescu2011state, wang2016solving, dvorkin2016optimizing,velay2019fully } is to produce a pre-contingency (or nominal) schedule for generators at minimal cost, such that it allows for feasible steady-state points of operation for a predefined set of \emph{credible contingencies}. A review of the SCOPF problem, its challenges and trends is available in \cite{capitanescu2011state}. 

The specification of the set of credible contingencies varies in academic works. Generally, a loss of up to one or two elements (generators and/or transmission lines) is considered. Interesting discussions about credible contingencies and reserve requirements can be found in \cite{bouffard2005umbrella,dvorkin2016optimizing} and the references therein. Security criteria and regulation for reserves also vary across independent system operators. A survey about the requirements for reserves across U.S. ancillary services can be found in \cite{zhou2016survey}. Without loss of generality, we consider the $N-1$ criterion for generators in this paper; that is, the system must operate under the loss of any individual generator.

Variants of the SCOPF problem include the \emph{corrective} case \cite{wang2016solving} where re-scheduling is possible, and the \emph{preventive} case where no re-dispatch occurs \cite{ li2008decomposed, dvorkin2016optimizing}. In this work, we consider preventive SCOPF with \emph{primary response} \cite{dvorkin2016optimizing}. In this framework, the synchronized generators must be able to automatically respond to contingencies
to restore the balance between loads and generation. 

Even though SCOPF is a nonlinear and nonconvex problem, for computational purposes, several authors adopt dc approximations \cite{bouffard2005umbrella, li2008decomposed, dvorkin2016optimizing,velay2019fully}. In such cases, stability constraints for the system 
 are generally expressed as power flow limits. In practical applications, the solution provided by the dc model can be checked for ac power flow feasibility. Then, iterative and/or heuristic procedures can be applied to further restrict the dc power flow constraints until feasibility is reached. 

The primary response of generators is explicitly modeled in \cite{restrepo2005unit} for a unit commitment application and in \cite{dvorkin2016optimizing, velay2019fully}, and \cite{karoui2010modeling} for SCOPF problems. In \cite{dvorkin2016optimizing, velay2019fully}, and \cite{restrepo2005unit} variables that represent the frequency drop in each contingency state were used to generate linear approximations of primary response. 
 These variables 
are multiplied by parameters that represent the sensitivity of generators to frequency changes. Such frequency regulation parameters are related to \emph{droop coefficients}. We refer the interested reader to the discussions about droop coefficients in \cite{dvorkin2016optimizing} (and references therein) where, particularly, the authors argue that the co-optimization of the droop coefficient and the SCOPF might save on costs. 
In this work, we have also opted for the dc power flow approximation with a linear model for primary response. In summary, for each contingency state, we have substituted the single variable representing frequency drop adopted in \cite{dvorkin2016optimizing} and \cite{velay2019fully} with a single global signal to generators (also a variable).


The SCOPF problem featuring automatic primary response of generators is a mixed-integer linear program (MILP) even under the dc relaxation. This is because the constraints for the automatic response of generators may lead to power outputs above generator limits \cite{restrepo2005unit}. To remedy this, we require binary variables for each generator and for each contingency state to determine whether a generator is producing according to the constraints for automatic response or at its limit. The size of the problem is generally large.  It is proportional to the number of contingencies since we are required to represent the network and the power flow variables for each post-contingency state.   

The Benders' decomposition approach, which has often been applied to solve energy planning problems \cite{li2008decomposed,wang2009contingency,bertsimas2013adaptive,wang2016solving}, is a natural candidate to tackle the preventive SCOPF problem. Generally, in Benders' approaches, the extensive formulation is recast into a master problem and subproblems. The master problem for power systems applications usually solves the nominal dispatch, and the subproblems represent the redispatch or corrective actions under contingencies and/or uncertain scenarios. An iterative procedure that involves solving the master problem and subproblems is performed. During this process, Benders' cuts for the violated subproblems are added into the master problem. The process continues until all subproblems are feasible. A valuable review on the Benders' decomposition method can be found in \cite{rahmaniani2017benders}. 

Unfortunately, preventive SCOPF imposes challenges for the application of traditional Benders' decomposition since the subproblems are nonconvex. The constraints that enforce the primary response of synchronized generators contain binary variables. Despite such challenges, a solution method inspired by \cite{bertsimas2013adaptive} considering nonconvex subproblems was provided in \cite{dvorkin2016optimizing}. However, the optimality for this method is not guaranteed. An alternative that ensures optimality is to recast the master problem to include the constraints for the primary response. This modification, however, does not scale well since the number of binary variables increases quadratically with the number of synchronized generators (assuming the $N-1$ security criteria for generators). 

 
In order to remedy the aforementioned limitations, we have devised an exact and scalable algorithm to tackle the preventive SCOPF problem with primary response. The focus of this work is on the computational and practical aspects of the solution methodology. The proposed decomposition scheme differs significantly from previously proposed solution methods. The outline of the method is as follows. 

In the master problem we consider a nominal optimal power flow problem that accounts for 
valid constraints for each contingency state. We initially disregard the network for contingency states and the disjunctions (binary variables) in the master problem. This approach alleviates the computational burden required for solving the master problem. During the iterative process only a small subset of the disjunctions and network constraints are introduced to the master problem by a column-and-constraint generation algorithm (CCGA) {\cite{BoZeng2011}}. In the proposed decomposition approach, the only optimization problem that is solved is the master problem. This is possible since we use: i) preprocessed structures based on the power transfer distribution factors (PTDF) that are useful both as feasibility checkers and as dedicated cuts in the post-contingency states, and ii) a numerical procedure that determines the post-contingency variables based solely on the nominal generation. The aforementioned preprocessed structures allow us to monitor the critical congested areas of the system. As it is necessary, these structures are transformed into constraints (that differ from Benders' cuts) that are added to the master problem. These cuts represent the network for the post-contingency states. Likewise, as it is necessary, we introduce the disjunctions (binary variables) representing the primary response model for a few contingency states into the master problem. We also propose a method to find 
high-quality primal solutions and a procedure that monitors the upper and lower bounds for the method.   

The rest of this paper is organized as follows: In Section II, the SCOPF model is introduced. The solution methodology is presented in Section III. Numerical experiments are reported in Section IV. Finally, this paper is concluded in Section V.

	
	
	
	
	  

\section{SCOPF with Primary response Formulation}

We assume a generic framework where a bid-based market for energy and reserve and/or unit commitment (UC) procedures have taken place hours before (or in the day before). 
We assume that, at the time the SCOPF is solved, the operator has precise forecasts for the few-hours-ahead non-dispatchable generation, and nodal net loads. For notational conciseness we assume that all generators are synchronized.

\subsection{Power Flow Constraints}

We use the following dc power flow constraints: 
\begin{align}
&\mathbf{A} \mathbf{f} + \mathbf{B} \mathbf{g} = \mathbf{{d}}                \label{eq.Master.EnergyBalance}\\
&\mathbf{f} = \mathbf{S}  \mathbf{\boldsymbol{\theta}}  \label{eq.sec.Kirc}\\ 
&-\mathbf{\overline{f}} \leq \mathbf{f} \leq \mathbf{\overline{f}} \label{eq.Master.PFLimit} \\
&\mathbf{\underline{g}}\leq\mathbf{g} \leq \mathbf{\overline{g}}\label{eq.Master.GenCap}\\[0.4em]
&\mathbf{A} \mathbf{f}_s + \mathbf{B} \mathbf{g}_s = \mathbf{{d}} & &\forall s \in \mathcal S \label{eq.Conting.EnergyBalance}\\
&\mathbf{f}_s = \mathbf{S}  \mathbf{\boldsymbol{\theta}}_s & &\forall s \in \mathcal S \label{eq.sec.Kirc.scenario}\\
&-\mathbf{\overline{f}} \leq \mathbf{f}_s \leq \mathbf{\overline{f}} \label{eq.PFLimit.scenario} & &\forall s \in \mathcal S\\
&\mathbf{g}_s \leq \mathbf{\overline{g}} & &\forall s \in \mathcal S.\label{eq.Scenario.GenCap}
\end{align}
Constraints \eqref{eq.Master.EnergyBalance}--\eqref{eq.Master.GenCap} model the power flow in the nominal state. Expression \eqref{eq.Master.EnergyBalance} represents nodal power balance under a dc power flow model, while Kirchhoff's second law is accounted for in (\ref{eq.sec.Kirc}). Transmission line limits and generator limits are enforced by (\ref{eq.Master.PFLimit}) and \eqref{eq.Master.GenCap}, respectively. In \eqref{eq.Master.GenCap}, we allow generation bounds $\mathbf{\underline{g}}$ and $\mathbf{\overline{g}}$ to be different from the minimum and maximum ($\mathbf{\hat{g}}$) set points for the generators due to commitment and/or operational constraints. Analogously to block \eqref{eq.Master.EnergyBalance}--\eqref{eq.Master.GenCap}, the set of constraints \eqref{eq.Conting.EnergyBalance}--\eqref{eq.Scenario.GenCap} model the power flow under each contingency state $s$. 

\subsection{Primary Response Model}
The primary response under contingency state $s$ is modeled by a global signal $n_s$ sent to all synchronized generators. This approach differs from those of \cite{dvorkin2016optimizing} and \cite{velay2019fully}, where variables representing frequency drops under contingency states are considered. We assume that the response of generator $i$ is proportional to its capacity $\hat{g}_i$ and also to a predefined coefficient $\gamma_i$ that is associated with the droop coefficient. Hence, under $s$, the automatic response of synchronized generator $i$ is given by $g_{s,i}-g_i = n_s \gamma_i \, \hat{g}_i$, with the additional constraints that $g_{s,i}\leq\overline{g}_i$.   Mathematically, we have
\begin{align}
&g_{s,i}= \min \{g_i + n_s \gamma_i \,\hat{g}_i,\;\overline{g}_i\} & &\forall i \in \mathcal{G}, \forall s \in \mathcal{S}, i\ne s
\label{disjunction_1}\\
&g_{s,s}=0 & & \forall s \in \mathcal{S} \label{n-1constraint}.
\end{align}
%
By using traditional MILP modeling techniques to rewrite
\eqref{disjunction_1}--\eqref{n-1constraint}, we obtain
\begin{align}
 &|g_{s,i}-g_i  - n_s \gamma_i \,\hat{g}_i| \leq\label{disj01}\\ & \qquad \qquad \qquad M(1-x_{s,i}) &\forall i \in \mathcal G, \forall s \in \mathcal S, i \ne s  \nonumber \\
 &g_i + n_s \gamma_i \,\hat{g}_i \geq \overline{g}_i (1-x_{s,i}) & \forall i \in \mathcal G, \forall s \in \mathcal S, i \ne s \\
 &g_{s,i} \geq \overline{g}_i (1-x_{s,i}) & \forall i \in \mathcal G, \forall s \in \mathcal S, i \ne s \\
 &n_s \in [0,1] & \forall s \in \mathcal S \\
 &x_{s,i} \in \{0,1\}& \forall i \in \mathcal G, \forall s \in \mathcal S \label{disj05}\\
 &g_{s,s}=0 &  \forall s \in \mathcal{S} \label{disj06}.
\end{align}
%

\subsection{Extensive Formulation for the SCOPF Problem}
\label{sec.EF}
The SCOPF problem is modeled as a MILP, where we minimize the cost of nominal generation in the objective function subject to constraints \eqref{eq.Master.EnergyBalance}--\eqref{eq.Scenario.GenCap} and \eqref{disj01}--\eqref{disj06}. For conciseness, we define $\mathcal{X} = [\mathbf{g},\mathbf{f},\mathbf{\boldsymbol{\theta}} ]$, $\mathcal{X}_s = [\mathbf{g}_s, \mathbf{f}_s, \mathbf{\boldsymbol{\theta}}_s  ]$, and $\mathcal{Y}_s = [\mathbf{g},\mathbf{g}_s, \mathbf{x}_s, n_s ] $. Let $\mathcal{Y}_s \in \mathcal{F}_s$ denote the disjunctions related to \eqref{disj01}--\eqref{disj06} while $\mathcal{X}\in\mathcal{E}$ and $\mathcal{X}_s\in\mathcal{E}_s$ denote the power flow constraints in nominal  (\eqref{eq.Master.EnergyBalance}--\eqref{eq.Master.GenCap}) and contingency states  (\eqref{eq.Conting.EnergyBalance}--\eqref{eq.Scenario.GenCap}) respectively. The extensive formulation for the SCOPF problem, labeled as EF method, is as follows
\begin{align}
\min_{\mathcal{X},[\mathcal{X}_{s},\mathbf{x}_{s}, n_s]_{s\in\mathcal{S}}}
\quad&  \mathbf{c}^{\top} \mathbf{g}     \label{obj.func}\\
\text{s.t.:}\quad &           \mathcal{X} \in \mathcal{E} \label{const.set01} \\
&            \mathcal{X}_s \in \mathcal{E}_s \label{const.set02} &\forall s \in \mathcal S \\
&            \mathcal{Y}_s \in \mathcal{F}_s \label{const.set03} &\forall s \in \mathcal S. 
\end{align}

\section{Solution Methodology}

In this work we focus on methods that guarantee optimality for the SCOPF problem. A Benders' decomposition with valid post-contingency constraints in the master problem is presented in Section \ref{Meth.Benders}. Preprocessed structures for feasibility checking and cut generation are presented in Section \ref{Meth.Lazy}, while a useful binary search is introduced in Section \ref{Meth.BinSearch}. The CCGA is described in Section \ref{Meth.CCGA}. A method for finding high-quality primal solutions 
is presented in Section \ref{Meth.Heur}.

\subsection{Modified Benders' Decomposition}
\label{Meth.Benders}

The intuitive Benders' decomposition approach for \eqref{obj.func}--\eqref{const.set03} is to define the master problem as the nominal schedule, associated with \eqref{obj.func}--\eqref{const.set01}, and the subproblems as the separable feasibility recourse problems enforcing \eqref{const.set02}--\eqref{const.set03} for each $s \in \mathcal{S}$. 
Unfortunately, this approach introduces nonconvexities to the subproblems \cite{dvorkin2016optimizing}, and thus, does not guarantee optimality. 

In order to ensure the convexity of the subproblems, and thus optimality for the method, we define the subproblems as feasibility-like problems for the constraints in \eqref{const.set02}
. As part of the modification, we also add the following valid post-contingency constraints to the master problem: 
\begin{align}
&  \mathbf{e}^{\top} \mathbf{g}_s   \label{valid.pstcont} = \mathbf{e}^{\top}\mathbf{d}  &\forall s \in \mathcal S.
\end{align}
The purpose of \eqref{valid.pstcont} is to strengthen the master problem with the necessary post-contingency condition that the total generation and the total load are equal. By also enforcing \eqref{eq.Scenario.GenCap} 
we guarantee that post-contingency generation is within bounds. We define the master problem as 
\begin{align}
\min_{\mathcal{X},[\mathbf{g}_{s},\mathbf{x}_{s}, n_s]_{s\in\mathcal{S}}}
\quad&  \mathbf{c}^{\top} \mathbf{g}   \label{master.benders} \\
\text{s.t.:}\quad &      \eqref{const.set01} \\
& \eqref{eq.Scenario.GenCap}, \eqref{const.set03},\eqref{valid.pstcont}& \forall s \in \mathcal{S}.\label{final.benders}
\end{align}

The subproblem for each $s\in\mathcal{S}$, where $\mathbf{g}_s^{(*)}$ is the solution determined in \eqref{master.benders}--\eqref{final.benders}, is then defined as 
\begin{align}
\min_{\mathbf{v}_s^+, \mathbf{v}_s^-, \mathbf{f}_{s}, \mathbf{g}_{s}, \mathbf{\boldsymbol{\theta}}}
\quad&  \mathbf{e}^{\top} ( \mathbf{v}_s^+ + \mathbf{v}_s^- )  \label{sub.benders01} \\ 
\text{s.t.:}\quad & \eqref{eq.sec.Kirc.scenario}, \eqref{eq.PFLimit.scenario}\\
& \mathbf{g}_s = \mathbf{g}_s^{(*)} & &:\boldsymbol{\mu}_s  \label{sub.benders03}\\
&\mathbf{A} \mathbf{f}_s + \mathbf{B} \mathbf{g}_s = \mathbf{{d}} + \mathbf{v}_s^+ - \mathbf{v}_s^-.& &  \label{sub.benders04}
\end{align}

A feasibility Benders' cut is then introduced to the master problem at each iteration, for each $s$ that is not feasible; that is, $\forall s \in \mathcal{S}$ such that $\mathbf{e}^{\top} ( \mathbf{v}^{+(*)}_s - \mathbf{v}^{-(*)}_s ) > \epsilon$, where $\epsilon$ is a tolerance level for the net load imbalance. The Benders' cut for $s$ is as follows:
%
$\mathbf{e}^{\top} ( \mathbf{v}^{+(*)}_s - \mathbf{v}^{-(*)}_s )   + \boldsymbol{\mu}^{\top} (\mathbf{g}_s - \mathbf{g}_s^{(*)}) \leq 0$. We label this approach as the BD method.

\subsection{Precomputation of Dedicated Cuts}
\label{Meth.Lazy}

In this subsection, an alternative method named BDDC is introduced. Unlike the BD method that involves subproblems that generate Benders' cuts, the BDDC method uses preprocessed structures as feasibility checkers and to generate cuts. These structures, that are also applied in the CCGA of Section \ref{Meth.CCGA}, are based on the PTDF formulation for dc power flow. 
 
 In the BDCC method, we have the same master problem \eqref{master.benders}--\eqref{final.benders} as the BD method. Thus, $\mathbf{g}_s^{(*)}$ is determined in \eqref{master.benders}--\eqref{final.benders}, where the primary response constraints \eqref{const.set03} and the post-contingency generation constraints \eqref{eq.Scenario.GenCap} 
 and \eqref{valid.pstcont} are enforced.
 
The aforementioned preprocessed structures are constructed directly from the PTDF-based formulation for the dc power flow. This formulation, for contingency state $s$, is as follows: 
\begin{align}
\min_{\mathbf{u}_s^-,\mathbf{u}_s^+} \quad  &\mathbf{0} \label{App2_01}\\
\text{s.t.} \quad  & -\mathbf{\overline{f}} + \mathbf{u}_s^- = \mathbf{K}_0 (\mathbf{{d}} - \mathbf{B} \mathbf{g}_s^{(*)} ) \label{preCut_01} = \mathbf{\overline{f}} - \mathbf{u}_s^+\\
 & \mathbf{u}_s^-,  \mathbf{u}_s^+ \geq \mathbf{0}. \label{eq.PTDFlast}
 \end{align}

In \eqref{preCut_01}, $\mathbf{g}_s^{(*)}$ is a solution determined by \eqref{master.benders}--\eqref{final.benders}  and $\mathbf{K}_0$ is the PTDF matrix. A similar description for \eqref{preCut_01} is presented in \cite{ardakani2013identification}. We highlight that \eqref{eq.Scenario.GenCap} 
and \eqref{valid.pstcont} are enforced in the master problem and therefore are not necessary in \eqref{App2_01}--\eqref{eq.PTDFlast}. As opposed to the subproblems of the BD method we do not allow nodal imbalance in \eqref{App2_01}--\eqref{eq.PTDFlast}. Thus, a given master problem solution $\mathbf{g}_s^{(*)}$ is feasible under contingency state $s$ if there is a feasible solution $\mathbf{u}_s^-$,  $\mathbf{u}_s^+$  for \eqref{App2_01}--\eqref{eq.PTDFlast}.

In this work we assume that we have more lines $|\mathcal{L}|$ than buses $|\mathcal{N}|$ and that there are no isolated buses. 
Under these assumptions and because we enforce \eqref{eq.Scenario.GenCap} 
and \eqref{valid.pstcont} in the master problem, we do not need to solve \eqref{App2_01}--\eqref{eq.PTDFlast} to check for feasibility in post-contingency states or to obtain cuts. Manipulating \eqref{preCut_01}, we derive the preprocessed structures:
\begin{align}
(\mathbf{u}_s^+):& \;\mathbf{K}_1 \mathbf{g}_s + \mathbf{k}_2\geq \mathbf{0} \label{Cut_01}\\
(\mathbf{u}_s^-):& \;\mathbf{K}_3 \mathbf{g}_s + \mathbf{k}_4\geq \mathbf{0} \label{Cut_02}.
\end{align}

Interestingly, the matrices $\mathbf{K}_1$, $\mathbf{K}_3$ and the vectors $\mathbf{k}_2$, $\mathbf{k}_4$ can be efficiently precomputed and are the same for all $s\in\mathcal{S}$. Another feature is that  \eqref{Cut_01} and \eqref{Cut_02} are directly associated with the transmission lines of the system
, relating either to the positive \eqref{Cut_01} or negative \eqref{Cut_02} limits. By inspecting solutions $\mathbf{g}_s$ on \eqref{Cut_01}--\eqref{Cut_02} it is possible to verify the existence of violated lines and the intensity of violations (in MW) under each contingency state. 


The algorithm for the BDDC method involves introducing rows of \eqref{Cut_01} and \eqref{Cut_02} as lazy constraints to problem \eqref{master.benders}--\eqref{final.benders}. We do not require that an optimal solution  $\mathbf{g}_{s}$ is found. Whenever a feasible (suboptimal) integer solution is determined by the solver, as a subroutine, we check the feasibility of $\mathbf{g}_{s}$  using \eqref{Cut_01} and \eqref{Cut_02}. We define $\alpha_{s,i}$ as the violation of line $l$, for contingency state $s$ and $\alpha$ as the largest violation among all transmission lines for all $s$. We then add to \eqref{master.benders}--\eqref{final.benders} as lazy constraints the rows of both \eqref{Cut_01} and \eqref{Cut_02} corresponding to violated lines such that $\alpha_{s,l}> \alpha/\beta_1$, where $\beta_1$ is a parameter. We stop the algorithm when $\alpha$ is smaller than a defined tolerance $\epsilon$. This procedure converges in finite steps since adding all rows of \eqref{Cut_01} and \eqref{Cut_02}, for every $s\in\mathcal{S}$, to \eqref{master.benders}--\eqref{final.benders} results in a problem which is equivalent to  \eqref{obj.func}--\eqref{const.set03}.

Unfortunately, the application of the BDDC method alone is not scalable since its master problem contains all the binary variables. Before proceeding to the proposed CCGA (Section \ref{Meth.CCGA}), we introduce next a useful binary search procedure.

\subsection{Numerical Procedure for Calculating $n_s$ and $\mathbf{g}_s$}
\label{Meth.BinSearch}

 Under the primary response model, the post-contingency generation for $s$; that is, $\mathbf{g}_s$, is defined by the combination of the nominal schedule $\mathbf{g}^{(*)}$ and the global signal $n_s^{(*)}$. Namely, given $\mathbf{g}^{(*)}$ and $n_s^{(*)}$, it is straightforward to compute $\mathbf{g}_s^{(*)}$ by applying the relations in  \eqref{disj01}--\eqref{disj06}. 
 
 Interestingly, $n_s^{(*)}$ can also be calculated from  $\mathbf{g}^{(*)}$. This is achieved by a binary search for $n_s$ for each $s\in\mathcal{S}$. The binary search is possible in this case since, for a fixed $\mathbf{g}^{(*)}$, each component of $\mathbf{g}_s$ is monotone with respect to $n_s$. Thus, despite the presence of the disjunctions, we only need to find the correct $n_s^{(*)}$ that results in a vector $\mathbf{g}_s^{(*)}$ that satisfies the total net demand. Given the fast convergence of the binary search, the tolerance $\varepsilon$ can be set to very small values even for large instances. This procedure 
 is as follows
 \begin{algorithm}[h!]
 	{\footnotesize \caption{{\small Binary Search($s$, $\varepsilon$, $\mathbf{g^{(*)}}$})} \label{Alg.Binary}
 		\begin{algorithmic}[1]
 			\State{Initialization: $j\leftarrow0$ and $n_s^{(0)}\leftarrow0.5$.}
 			\For{$\forall i\in\mathcal{G}, i\ne s$}
 			\If{$g_i^{(*)} + n_s^{(j)} \gamma_i \,\hat{g}_i \geq \overline{g}_i$ }\textbf{:} $g_{s,i}^{(j)} \leftarrow \overline{g}_i$
 			\Else{}\textbf{:} $g_{s,i}^{(j)} \leftarrow g_i^{(*)} + n_s^{(j)} \gamma_i \,\hat{g}_i$
 			\EndIf.
 			\EndFor
 			\State{$g_{s,s}^{(j)}\leftarrow0$;  $e_s \leftarrow (\mathbf{e}^{\top} \mathbf{g}_s  - \mathbf{e}^{\top}\mathbf{d})$.}
 			
 			\If{$|e_s| \le \varepsilon$}\textbf{:} $n_s^{(*)} \leftarrow n_s^{(j)}$,  $g_{s,i}^{(*)} \leftarrow g_{s,i}^{(j)}, \forall i\in\mathcal{G}$ and \textbf{BREAK} 
 			 
 			\ElsIf{$e_s>0$ }\textbf{:} $n_s^{(j+1)}\leftarrow n_s^{(j)}/2$
 			\Else \textbf{:} $n_s^{(j+1)}\leftarrow (1+n_s^{(j)})/2$
 			\EndIf.
 			\State{$j\leftarrow j+1$; Go to step 2.}
 		\end{algorithmic}
 	} 
 \end{algorithm}

 \subsection{Column-and-Constraint-Generation Algorithm }
 \label{Meth.CCGA}
 
We define the master problem for the CCGA as follows 
\begin{align}
z = \min_{\mathcal{X},[\mathbf{g}_{s}^{'}]_{s\in\mathcal{S}},[\mathbf{x}_{s},n_s]_{s\in\mathbb{S}}}
\quad&   \mathbf{c}^{\top} \mathbf{g}   \label{master.ccga} \\
\text{s.t.:}\quad &      \eqref{const.set01} \\
& \mathbf{g}_{s}^{'} -\mathbf{g} \leq  \mathbf{\overline{r}} & \forall s \in \mathcal S \label{RampingCCGA}\\
& \eqref{eq.Scenario.GenCap}, \eqref{disj06}, \eqref{valid.pstcont}& \forall s \in \mathcal{S}\label{semifinal.ccga}\\
& \eqref{const.set03}& \forall s \in \mathbb{S}.\label{final.ccga}
\end{align}

Note that, as opposed to the BDDC method, in \eqref{final.ccga} we define a different set of contingency states $\mathbb{S}$ for the disjunctive constraints (starting with $\mathbb{S}=\emptyset$). We also abuse notation by using $\mathbf{g}_{s}^{'}$ as a provisional variable for the post-contingency generation replacing $\mathbf{g}_{s}$ in \eqref{RampingCCGA} and \eqref{semifinal.ccga}. We performed this substitution to make explicit that $\mathbf{g}_{s}$ is not determined in \eqref{master.ccga}--\eqref{final.ccga} for the entire iterative process. The determination of $\mathbf{g}_{s}$ in \eqref{master.ccga}--\eqref{final.ccga} would only be possible in the presence of the disjunctive constraints \eqref{const.set03} for $s$. These disjunctive constraints are not initially present in \eqref{master.ccga}--\eqref{final.ccga} for computational purposes. In fact, the determination of $\mathbf{g}_{s}$ is performed by Algorithm \ref{Alg.Binary}, which requires only $\mathbf{g}$ as an input. 
The purpose of $\mathbf{g}_{s}^{'}$ in \eqref{RampingCCGA} and \eqref{semifinal.ccga} is to guarantee that $\mathbf{g}$ is determined in such a way that Algorithm \ref{Alg.Binary} is capable of enforcing the primary response compatibility to $\mathbf{g}_{s}$, while meeting the global demand. That is, for each $s$, $|\mathbf{e}^{\top} \mathbf{g}_s  - \mathbf{e}^{\top}\mathbf{d}|\leq\epsilon$. 

In order to verify the above claim, note that by \eqref{RampingCCGA} and \eqref{eq.Scenario.GenCap}, $g_{s,i}^{'}\leq \min{\{\overline{g}_i, g_i + \gamma_i \hat{g}_i\}}$ for each $i$ and $s$, where $g_{s,i}^{'}$ is the $i$-th element of $\mathbf{g}_{s}^{'}$. Defining $n_s=0$ in Algorithm \ref{Alg.Binary} implies $\mathbf{g}_{s}=\mathbf{g}$, except for $g_{s,s}=0$. If instead we set $n_s=1$ then $g_{s,i} = \min{\{\overline{g}_i, g_i + \gamma_i \hat{g}_i\}} \geq g_{s,i}^{'}$ for each $i$ and $s$, with $i\ne s$. For $i=s$, we have that $g_{s,s}^{'}=g_{s,s}=0$. Thus, since $\mathbf{g}_{s}^{'}$ meets the global demand, it is always the case that $\mathbf{e}^{\top} \mathbf{g}_s  \geq \mathbf{e}^{\top} \mathbf{g}_s^{'} = \mathbf{e}^{\top}\mathbf{d}$ by choosing $n_s=1$. By the monotonicity and continuity of $g_{s,i}$ with respect to $n_s$ for a given $g_i$, there is a value $n_s^{(*)}$ that results in $\mathbf{g}_{s}^{(*)}$ that satisfies the global demand and preserves the primary response model.


At each iteration $j$ of the CCGA, we solve the master problem \eqref{master.ccga}--\eqref{final.ccga} to obtain $\mathbf{g}^{(j)}$ and $z^{(j)}$. Then, for each $s\in\mathcal{S}$, we perform the binary search (Algorithm \ref{Alg.Binary}) to define $\mathbf{g}_{s}^{(j)}$ according to the primary response model. Next, for all $s$, we check feasibility of the solutions $\mathbf{g}_{s}^{(j)}$  using \eqref{Cut_01} and \eqref{Cut_02}. We use $\alpha_{s,i}^{(j)}$ to define the violation of each line $l$, for each contingency state $s$ and we use $\alpha^{(j)}$ as the largest violation among all transmission lines for all $s$. We identify the contingency state $s^{(j)}$ that contains the most violated line. If $s^{(j)}\in \mathbb{S}$ we skip the rest of this step. Otherwise we set $\mathbb{S}=\mathbb{S}\cup s^{(j)}$ which means including the disjunctions \eqref{const.set03} for $s^{(j)}$ into \eqref{master.ccga}--\eqref{final.ccga}. 

We also introduce to \eqref{master.ccga}--\eqref{final.ccga} the rows of both \eqref{Cut_01} and \eqref{Cut_02} corresponding to violated lines using two criteria: i) For the post-contingency states $s\in \mathbb{S}$ we include the lines where $\alpha^{(j)}_{s,l}> \alpha^{(j)}/\beta_1$. ii) For the post-contingency states $s\notin \mathbb{S}$ we include the lines where $\alpha^{(j)}_{s,l}> \alpha^{(j)}/\beta_2$. The objective of this step is to enforce the network constraints for critical lines in post-contingency states. Typically, $\beta_1 > \beta_2$. We are stricter with the states $s\in \mathbb{S}$ since defining tight parameters for contingency states without corresponding disjunctions may lead to the inclusion of many constraints at a time. A user defined tolerance $\epsilon $ (in MW) for maximum line violation is used to stop the iterative process. The algorithm is as follows:
\begin{algorithm}[H]
	{\footnotesize \caption{{\small CCGA($\varepsilon$, $\epsilon$, $\beta_1$, $\beta_2$)} }\label{Alg.CCGA}
		\begin{algorithmic}[1]
			\State{Initialization: $j\leftarrow0$, $\mathbb{S} \leftarrow \emptyset$.}
			\State{Solve: \eqref{master.ccga}--\eqref{final.ccga} to obtain $\mathbf{g}^{(j)}$}
			\For{$\forall s\in\mathcal{S}$}
			\State{Use Algorithm \ref{Alg.Binary} to obtain: $z^{(j)}$, $\mathbf{g}^{(j)}_s, n^{(j)}_s$}
				\For{$\forall l\in\mathcal{L}$}
					\State{Compute $\alpha^{(j)}_{s,l}$ as the maximum violation of \eqref{Cut_01} and \eqref{Cut_02}}
				\EndFor
			\EndFor
			\State{Compute $\alpha^{(j)}$ and identify the state $s^{(j)}$ of the most violated line}
			\If{$s^{(j)}\notin \mathbb{S}$}	set $\mathbb{S}=\mathbb{S}\cup s^{(j)}$ (add \eqref{const.set03} for $s^{(j)}$ to \eqref{master.ccga}--\eqref{final.ccga}) 
			\EndIf
			\For{$\forall s\in\mathcal{S}$}
				\If{$s\in\mathbb{S}$} $\beta\leftarrow\beta_1$
				\Else { $\beta\leftarrow\beta_2$}
				\EndIf
				\For {$\forall l\in\mathcal{L}$}
					\If{$\alpha^{(j)}_{s,l}> \alpha^{(j)} / \beta$ \textbf{AND if} cuts \eqref{Cut_01}--\eqref{Cut_02} for the pair ($l$,$s$) are not yet included} \State{Add \eqref{Cut_01}--\eqref{Cut_02} for contingency state $s$, line $l$}
					\EndIf
				\EndFor
			\EndFor
			\If{$\alpha^{(j)}<\epsilon$} $z^{(*)}\leftarrow z^{(j)}$, $\mathbf{g}^{(*)}\leftarrow \mathbf{g}^{(j)}$; \textbf{BREAK}
			\Else { $j\leftarrow j+1$; Go to step 2.}
			\EndIf
		\end{algorithmic}
	} 
\end{algorithm}

\subsection{Finding High-Quality Primal Solutions and Monitoring the Optimality Gap using the CCGA}
\label{Meth.Heur}
Because very large cases might still impose computational challenges, we propose a procedure for finding feasible primal solutions. This procedure restricts the disjunctions in \eqref{const.set03} to a subset of synchronized generators $\mathcal{H}\subset\mathcal{G}$. The generators in $\mathcal{G}\setminus\mathcal{H}$ respond with $g_{s,i}-g_i = n_s \gamma_i \, \hat{g}_i$; that is, we define $x_{s,i} = 1, \forall s\in\mathcal{S}, \forall i\in\mathcal{G}\setminus\mathcal{H}$, $i\ne s$. 

The following criterion is used for defining $\mathcal{H}$. We rank the synchronized generators according to a “cost/limit” index ($c_i/\overline{g}_i$) and define $\mathcal{H}$ as the $p\%$ generators with lowest ranks, where $p$ is a parameter. The objective value of the problem using this approach is denoted as  $z^{(*)}_{p}$. Note that $z^{(*)}_{100}=z^{(*)}$.

This approach reduces the number of binary variables and thus the complexity of the problem. It is then a tool for finding upper bounds for \eqref{obj.func}--\eqref{const.set03}. If we apply Algorithm \ref{Alg.CCGA} using the proposed primal method; i.e., setting $x_{s,i} = 1, \forall s\in\mathcal{S}, \forall i\in\mathcal{G}\setminus\mathcal{H}$ in \eqref{master.ccga}--\eqref{final.ccga} we obtain $ub = z^{(*)}_{p}$ as a valid upper bound. If the problem is infeasible for $p$ then $z^{(*)}_{p}\leftarrow\inf$.


Note, however, that a procedure that monitors the optimality gap is still required. A lower bound for \eqref{obj.func}--\eqref{const.set03} is not obtained for $p<100$. Conversely, solving for $p=100$ generates upper bounds only after a feasible solution is found. This typically occurs in the later iterations when the tolerance $\epsilon$ for all lines in every contingency state is achieved. 

We propose a simple strategy that monitors the bounds. Note that solving the SCOPF with the CCGA for low values of $p$ tends to be faster than for high values of $p$. Thus, we start $p=0$ and increase it sequentially. The solution for each $p<100$ provides 
an upper bound for \eqref{obj.func}--\eqref{const.set03}. As a parallel procedure, we solve for $p=100$ to obtain valid lower bounds. Namely, for each iteration $j$ of Algorithm \ref{Alg.CCGA} for $p=100$, a valid lower bound is defined as the best bound provided by the solver. This procedure monitors the gap efficiently.

%
%
%

\section{Computational Experiments}
We compared the proposed CCGA with two solution methods: EF and BDDC described in Sections \ref{sec.EF}  and \ref{Meth.Lazy}.

We performed simulations for various values of $\boldsymbol{\gamma}$, $\beta_1$, and $\beta_2$. Our results indicate that varying the parameters may impact the performance of the CCGA. However, the dominance of the CCGA over other methods (EF and BDDC) was a constant, despite the parameterization. We have reported results for $\beta_1=5$, $\beta_2=1.2$, and $\gamma_i=0.05$ for all $i\in\mathcal{G}$.



We used Gurobi 8.1.1 under the modeling package JuMP 0.18.5 for Julia Language 0.6.4 on a Xeon E5-2680 processor at 2.5 GHz and 128 GB of RAM. We set the optimality gap of Gurobi to $0.5$\% for the EF method and BDDC method as well as for the master problem of the CCGA. The maximum line violation was set to $\epsilon=0.05$ MW and the precision of Algorithm \ref{Alg.Binary} to $\varepsilon=10^{-10}$ MW.%

The data are based on modified versions of the benchmark systems presented in \cite{PGLIB}. 
The size of the instances for the EF method, after Gurobi's presolve, are reported in Table \ref{tab:Inst.size}.
%
%
%
\begin{table}[t!]
	\centering\scriptsize
	\caption{Instance Size for the EF method after Presolve}
	\begin{tabular}{lr|r|r}
		\toprule
		\noalign{\vskip 0.35mm}
	    \multirow{2}{0.40cm}{System}	& \multicolumn{1}{c|}{Continuous} &  \multicolumn{1}{c|}{Binary} &  \multicolumn{1}{c}{Linear}  \\
										&  \multicolumn{1}{c|}{Variables} & \multicolumn{1}{c|}{Variables}	&  \multicolumn{1}{c}{Constraints}    \\
		\midrule
		118 IEEE	     &10,604		&2,862		&  19,137	\\
		1354 PEGASE      &323,571 		&63,455	   	&  513,677\\
		1888 RTE         &387,979  		&79,032   	&  624,780 \\
		1951 RTE         & 563,273  	&149,370   	&  1,010,994\\
		2848 RTE         & 1,104,192   	&276,822   	&  1,934,115 	\\
		2868 RTE         & 1,284,568 	&348,036    &  2,328,081\\
		6468 RTE         & 5,067,009  	&1,563,640 	&  9,756,668\\
		\bottomrule\\
			\vspace{-0.47cm}
	\end{tabular}%
	\label{tab:Inst.size}%
\end{table}

\subsection{Solution Method Comparison}

 Table \ref{tab:CPU_TIME} provides the computational times for selected methods and the number of iterations for the CCGA method in parentheses. The CCGA dominates the other methods, which were only able to solve the 118 IEEE case within a reasonable time limit. For this instance, the CCGA required less than one third of the time of the BDDC method and less than one fifth of the time of the EF mehod. The only other instance that the BDDC method was able to solve in less than 4 hours was the 1888 RTE. The CCGA was more than 200 times faster for this instance. Interestingly, the number of iterations required by the CCGA is generally small, implying that the CCGA solved far less complex MILPs than the other methods. The only instance that posed difficulties for the CCGA was the 6468 RTE that contains more than 6,000 buses, 1,200 generators, and 9,000 transmission lines. Nevertheless, a solution for the optimality gap of $0.5\%$ was achieved in less than 3 hours. 

Interestingly, as reported next, it is possible to determine high-quality solutions in competitive computational times for large systems by applying the primal method of Section \ref{Meth.Heur}.

\begin{table}[t]
	\centering\scriptsize
	\caption{Comparative CPU times (s) and number of iterations}
	\begin{tabular}{lr|r|r}
		\toprule
		\multirow{3}{0.40cm}{System} & \multicolumn{3}{c}{Solution Method}\\
		\noalign{\vskip 0.35mm}
		\cline{2-4} \noalign{\vskip 0.95mm}
		& EF & BDDC & CCGA  \\
		\midrule
		118 IEEE	     &46.2		&27.6		&  8.0$\;\,\,\,$(4)	\\
		1354 PEGASE	     &T 		&T	   		&  138.0 (12)\\
		1888 RTE         &T  		&2899.0   	&  14.0$\;\,\,\,$(4) \\
		1951 RTE         & T  		&T   	 	&   16.2$\;\,\,\,$(4)\\
		2848 RTE         & T   		&T   	 	&   26.6$\;\,\,\,$(4)\\
		2868 RTE         & T 		&T      	&   31.0$\;\,\,\,$(3)\\
		6468 RTE         & T  		&T  	 	&   7881.6 (18)\\
		\bottomrule\\
		\vspace{-0.47cm}\\
		\multicolumn{4}{l}{T - Time limit of 4 hours exceeded.}
	\end{tabular}%
	\label{tab:CPU_TIME}%
\end{table}
\vspace{-0.18cm}
\subsection{Finding Primal Solutions and Determining Bounds}

The method of Section \ref{Meth.Heur} for defining primal solutions was applied for the 6468 RTE and 1354 PEGASE systems. CCGA was used to solve the SCOPF problem for different values of $p$ to an optimality gap of $0.5\%$. The results are summarized in Tables \ref{tab:Heur1} and \ref{tab:Heur2}. Columns 1 and 2 present the cost and the relative cost gap for each $p$ with respect to the cost achieved by $p=100$. Columns 3 and 4 report the required computational time and number of iterations.
\begin{table}[t!]
	\centering\scriptsize
	\caption{Primal Approach for the 6468 RTE System}
	\begin{tabular}{rr|c|r|r   }
		\toprule
		\noalign{\vskip 0.35mm}
		\multirow{2}{0.21cm}{$p \;\;\quad$}& \multicolumn{1}{c|}{Cost} & \multicolumn{1}{c|}{Cost Gap} &  \multicolumn{1}{c|}{Time} &  \multicolumn{1}{c}{Iter.}   \\
		& \multicolumn{1}{c|}{($10^3$ \$)}& \multicolumn{1}{c|}{(\%)} &  \multicolumn{1}{c|}{(s)} &  \multicolumn{1}{c}{(\#)}   \\
		\midrule
		100	    &   1624.8    &0.00			&7881.6		&  18\\
		10      &   1625.8    &0.02  		&813.5   	 	&   17   \\	
		0       &   1628.7    &0.25			&481.2     	&   16   \\
		\bottomrule\\
		\vspace{-0.47cm}\\
	\end{tabular}%
	\label{tab:Heur1}%
\end{table}
\begin{table}[t!]
	\centering\scriptsize
	\caption{Primal Approach for the 1354 PEGASE System}
	\begin{tabular}{rr|c|r|r   }
		\toprule
		\noalign{\vskip 0.35mm}
			\multirow{2}{0.21cm}{$p \;\;\quad$}& \multicolumn{1}{c|}{Cost} & \multicolumn{1}{c|}{Cost Gap} &  \multicolumn{1}{c|}{Time} &  \multicolumn{1}{c}{Iter.}   \\
		& \multicolumn{1}{c|}{($10^3$ \$)}& \multicolumn{1}{c|}{(\%)} &  \multicolumn{1}{c|}{(s)} &  \multicolumn{1}{c}{(\#)}   \\
		\midrule
		100	    &   1190.6    &0.00			&138.0		&  12\\
		50      &   1192.6    &0.17  		&64.7 		&  11  \\
		10      &   1195.4    &0.40  		&21.1   	 	&   11   \\	
		0       &   1208.3    &1.49			&9.4     	&   10   \\
		\bottomrule\\
		\vspace{-0.47cm}\\
	\end{tabular}%
	\label{tab:Heur2}%
\end{table}

For the 6468 RTE system (Table \ref{tab:Heur1}) the result is quasi-optimal even for $p=0$. For the 1354 PEGASE system (Table \ref{tab:Heur2}) the solution for $p=0$ is already competitive, and required $9.4$ seconds only. By increasing the complexity of the problem to $p=10$, the cost gap is reduced by more than $1\%$ for a reasonable solution time of $21.1$ seconds. For $p=50$, the CCGA required $64.7$ seconds to converge, achieving a negligible cost gap of $0.17$. 

\begin{figure}[!t]
	\centering
	\includegraphics{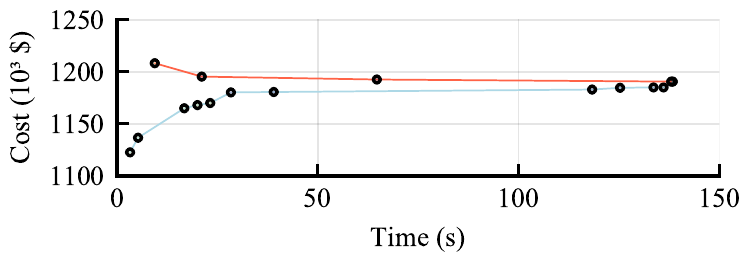}
	\caption{Bounds for the 1354 PEGASE system.}
	\label{fig_sim}
\end{figure}

Despite the good results for small values of $p$, the cost gap is not observable before solving for $p=100$. Thus, 
 we adopted the strategy of Section \ref{Meth.Heur} for obtaining bounds. We used multi-threading to solve problems in parallel. In the first thread we solved a sequence of problems for increasing values of $p$, starting with $p=0$. We have stored the costs and times for the solutions of each $p$. In the second thread we solved for $p=100$ and recorded solving time and the best bound of each iteration provided by Gurobi. A convergence plot from applying this method to the 1354 PEGASE system is illustrated in Fig. \ref{fig_sim}.  
 The proposed strategy yields the true optimality gap and is a useful decision-making tool for system operators. 




\section{Conclusion}
We presented an exact and scalable column-and-constraint-generation algorithm for the SCOPF problem with primary response of generators. Under the proposed framework, we add the disjunctions as necessary in an iterative process that does not involve subproblems. This is possible by a scheme that involves a master problem with valid post-contingency constraints, preprocessed structures that serve both as feasibility checkers and delayed cuts, and a numerical procedure that reduces the complexity of the master problem by exogenously calculating the nonconvex primary response.  We also proposed a procedure for finding high-quality primal solutions that helps monitor the bounds for the method. As shown by the computational experiments, this approach scales to large instances of the SCOPF problem with primary response. 

\bibliographystyle{IEEEtran}
\bibliography{References.bib}

\begin{thebibliography}{10}
\providecommand{\url}[1]{#1}
\csname url@samestyle\endcsname
\providecommand{\newblock}{\relax}
\providecommand{\bibinfo}[2]{#2}
\providecommand{\BIBentrySTDinterwordspacing}{\spaceskip=0pt\relax}
\providecommand{\BIBentryALTinterwordstretchfactor}{4}
\providecommand{\BIBentryALTinterwordspacing}{\spaceskip=\fontdimen2\font plus
\BIBentryALTinterwordstretchfactor\fontdimen3\font minus
  \fontdimen4\font\relax}
\providecommand{\BIBforeignlanguage}[2]{{%
\expandafter\ifx\csname l@#1\endcsname\relax
\typeout{** WARNING: IEEEtran.bst: No hyphenation pattern has been}%
\typeout{** loaded for the language `#1'. Using the pattern for}%
\typeout{** the default language instead.}%
\else
\language=\csname l@#1\endcsname
\fi
#2}}
\providecommand{\BIBdecl}{\relax}
\BIBdecl

\bibitem{alsac1974optimal}
O.~Alsac and B.~Stott, ``Optimal load flow with steady-state security,''
  \emph{IEEE Trans. Power App. Syst.}, vol. PAS-93, no.~3, pp. 745--751, May
  1974.

\bibitem{bouffard2005umbrella}
F.~Bouffard, F.~D. Galiana, and J.~M. Arroyo, ``Umbrella contingencies in
  security-constrained optimal power flow,'' in \emph{15th Power systems
  computation conference, PSCC}, vol.~5, 2005.

\bibitem{li2008decomposed}
Y.~Li and J.~D. McCalley, ``Decomposed {SCOPF} for improving efficiency,''
  \emph{IEEE Trans. Power Syst.}, vol.~24, no.~1, pp. 494--495, Feb. 2009.

\bibitem{capitanescu2011state}
F.~Capitanescu, J.~M. Ramos, P.~Panciatici, D.~Kirschen, A.~M. Marcolini,
  L.~Platbrood, and L.~Wehenkel, ``State-of-the-art, challenges, and future
  trends in security constrained optimal power flow,'' \emph{Elect. Power Syst.
  Res.}, vol.~81, no.~8, pp. 1731--1741, Aug. 2011.

\bibitem{wang2016solving}
Q.~Wang, J.~D. McCalley, T.~Zheng, and E.~Litvinov, ``Solving corrective
  risk-based security-constrained optimal power flow with lagrangian relaxation
  and benders decomposition,'' \emph{Int. J. Elec. Power}, vol.~75, pp.
  255--264, Feb. 2016.

\bibitem{dvorkin2016optimizing}
Y.~Dvorkin, P.~Henneaux, D.~S. Kirschen, and H.~Pand{\v{z}}i{\'c}, ``Optimizing
  primary response in preventive security-constrained optimal power flow,''
  \emph{IEEE Syst. Journal}, vol.~12, no.~1, pp. 414--423, Mar. 2016.

\bibitem{velay2019fully}
M.~Velay, M.~Vinyals, Y.~Besanger, and N.~Retiere, ``Fully distributed security
  constrained optimal power flow with primary frequency control,'' \emph{Int.
  J. Elec. Power}, vol. 110, pp. 536--547, Sep. 2019.

\bibitem{zhou2016survey}
Z.~Zhou, T.~Levin, and G.~Conzelmann, ``Survey of {U.S.} ancillary services
  markets,'' Argonne National Lab.(ANL), Argonne, IL (United States), Tech.
  Rep., 2016.

\bibitem{restrepo2005unit}
J.~F. Restrepo and F.~D. Galiana, ``Unit commitment with primary frequency
  regulation constraints,'' \emph{IEEE Trans. Power Syst.}, vol.~20, no.~4, pp.
  1836--1842, Nov. 2005.

\bibitem{karoui2010modeling}
K.~Karoui, H.~Crisciu, and L.~Platbrood, ``Modeling the primary reserve
  allocation in preventive and curative security constrained {OPF},'' in
  \emph{Proc. IEEE PES Trans. Distrb. Conf. Expo}.\hskip 1em plus 0.5em minus
  0.4em\relax IEEE, Apr., 2010, pp. 1--6.

\bibitem{wang2009contingency}
J.~Wang, M.~Shahidehpour, and Z.~Li, ``Contingency-constrained reserve
  requirements in joint energy and ancillary services auction,'' \emph{Trans.
  Power Syst.}, vol.~24, no.~3, pp. 1457--1468, Aug. 2009.

\bibitem{bertsimas2013adaptive}
D.~Bertsimas, E.~Litvinov, X.~A. Sun, J.~Zhao, and T.~Zheng, ``Adaptive robust
  optimization for the security constrained unit commitment problem,''
  \emph{Trans. Power Syst.}, vol.~28, no.~1, pp. 52--63, Feb. 2013.

\bibitem{rahmaniani2017benders}
R.~Rahmaniani, T.~G. Crainic, M.~Gendreau, and W.~Rei, ``The {Benders}
  decomposition algorithm: {A} literature review,'' \emph{Eur. J. Oper. Res.},
  vol. 259, no.~3, pp. 801--817, Jun. 2017.

\bibitem{BoZeng2011}
B.~Zeng and L.~Zhao, ``Solving two-stage robust optimization problems using a
  column-and-constraint generation method,'' \emph{Oper. Res. Lett.}, vol.~41,
  no.~5, pp. 457--461, Sep. 2013.

\bibitem{ardakani2013identification}
A.~J. Ardakani and F.~Bouffard, ``Identification of umbrella constraints in
  dc-based security-constrained optimal power flow,'' \emph{IEEE Trans. Power
  Syst.}, vol.~28, no.~4, pp. 3924--3934, Nov. 2013.

\bibitem{PGLIB}
S.~Babaeinejadsarookolaee \emph{et~al.}, ``The power grid library for
  benchmarking ac optimal power flow algorithms,'' \emph{arXiv preprint
  arXiv:1908.02788}, 2019.

\end{thebibliography}

\end{document}